\documentclass[11pt]{article}
\usepackage{amssymb,amsfonts,amsmath,latexsym,epsf,tikz,url}

\newtheorem{theorem}{Theorem}[section]
\newtheorem{proposition}[theorem]{Proposition}

\newtheorem{definition}[theorem]{Definition}

\newcommand{\proof}{\noindent{\bf Proof.\ }}
\newcommand{\qed}{\hfill $\square$\medskip}

\begin{document}

\title{Characterization of graphs with distinguishing number equal list distinguishing number}

\author{
Saeid Alikhani  $^{}$\footnote{Corresponding author}
\and
Samaneh Soltani
}

\date{\today}

\maketitle

\begin{center}
Department of Mathematics, Yazd University, 89195-741, Yazd, Iran\\
{\tt alikhani@yazd.ac.ir, s.soltani1979@gmail.com}
\end{center}

\begin{abstract}
The distinguishing number $D(G)$ of a graph $G$ is the least integer $d$
such that $G$ has an vertex labeling   with $d$ labels  that is preserved only by a trivial automorphism. A list assignment to $G$ is an
assignment $L = \{L(v)\}_{v\in V (G)}$ of lists of labels to the vertices of $G$. A distinguishing $L$-labeling  of $G$ is a distinguishing labeling of $G$ where the label of each vertex $v$ comes
from $L(v)$. The list distinguishing number of $G$, $D_l(G)$ is the minimum $k$ such that every list assignment to $G$ in which $|L(v)| = k$ for all $v \in V (G)$ yields a distinguishing $L$-labeling
of $G$. In this paper, we determine the
list-distinguishing number for two families of graphs. We also characterize graphs with the distinguishing number equal the list distinguishing number. Finally, we show that this characterization works for other list numbers of a graph.   
\end{abstract}

\noindent{\bf Keywords:}   Distinguishing number;   list-distinguishing labeling; list distinguishing chromatic number.

\medskip
\noindent{\bf AMS Subj.\ Class.:} 05C15, 05E18

\section{Introduction}

Let $G = (V ,E)$ be a simple  graph.  We use the standard graph notation (\cite{Sandi}).  
The set of all \textit{automorphisms} of $G$, with the operation of composition of permutations, is a permutation group
on $V$ and is denoted by ${\rm Aut}(G)$.  
A labeling of $G$, $\phi : V \rightarrow \{1, 2, \ldots , r\}$, is  \textit{$r$-distinguishing}, 
if no non-trivial  automorphism of $G$ preserves all of the vertex labels.
In other words,  $\phi$ is $r$-distinguishing if for every non-trivial $\sigma \in {\rm Aut}(G)$, there
exists $x$ in $V$ such that $\phi(x) \neq \phi(\sigma(x))$. 
The distinguishing number of a graph $G$ has been defined in  \cite{Albert} and  is the minimum number $r$ such that $G$ has a labeling that is $r$-distinguishing.  
The introduction of the distinguishing number  was a great success; by now about one hundred papers were written motivated by this seminal paper! The core of the research has been done on the invariant $D$ itself, either on finite \cite{Chan, immel-2017, Kim}  or infinite graphs~\cite{estaji-2017, lehner-2016, smith-2014}; see also the references therein.

Recently Ferrara et al. \cite{Ferrara} extended the notion of a distinguishing labeling to a list distinguishing labeling.
 A \textit{list assignment} to $G$ is an
assignment $L = \{L(v)\}_{v\in V (G)}$ of lists of labels to the vertices of $G$. A \textit{distinguishing $L$-labeling}  of $G$ is a distinguishing labeling of $G$ where the label of each vertex $v$ comes
from $L(v)$. The \textit{list distinguishing number}  of $G$, $D_l(G)$ is the minimum $k$ such that every list assignment to $G$ in which $|L(v)| = k$ for all $v \in V (G)$ yields a distinguishing $L$-labeling
of $G$. Since all of the lists can be identical, we observe that $D(G) \leq D_l(G)$.
In some cases, it is easy to show that the list-distinguishing number can  equal the distinguishing number. For example, it is not difficult to see that $D(K_n) = n = D_l(K_n)$, 
$D(K_{n,n}) = n + 1 = D_l(K_{n,n})$ and $D_l(C_n) = D(C_n) = 2$ (\cite{Ferrara}). In particular  Ferrara et al.  \cite{Ferrara2} extend an enumerative technique of Cheng \cite{cheng}, to show that for any tree $T$, $D_l(T ) = D(T )$. Ferrara et al. \cite{Ferrara} asked the
following question at the end of their paper.

\medskip
\noindent \textbf{Question} Does there exist a graph $G$ such that $D(G) \neq D_l(G)$?

\medskip
Amusingly, Ferrara feels that no such graph $G$ exists, while Gethner
believes this Question  can be answered in the affirmative. 

\medskip
In this paper we first compute the
list-distinguishing number   for friendship and book graphs. We  also state a necessary and sufficient condition for  graph $G$ satisfying   $D_l(G ) = D(G)$ in Section 3. 
Finally in Section 4, we extend the results of Section 3 for other  list numbers of a graph. 

\section{List distinguishing number of friendship and book graphs}

\begin{figure}[ht]
\begin{center}
		\hspace{.7cm}
		\includegraphics[width=0.5\textwidth]{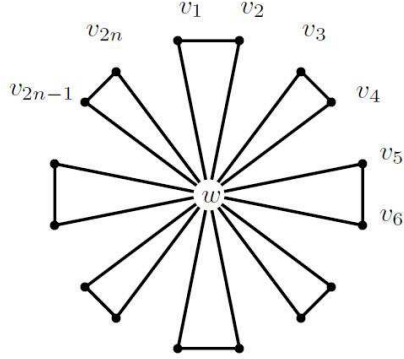}
	\caption{\label{friend} Friendship graph $F_n$.}
\end{center}
\end{figure}

 In this section, we consider the friendship graphs and the book graphs and compute their  list-distinguishing number. We begin with the friendship graph. 
The friendship graph $F_n$ $(n\geqslant 2)$ can be constructed by joining $n$ copies of the cycle graph $C_3$ with a common vertex (see Figure \ref{friend}). 

\begin{theorem}{\rm \cite{Alikhani}}
For every  $n\geq 2$, $D(F_n)= \lceil \dfrac{1+\sqrt{8n+1}}{2}\rceil.$  
\end{theorem}
	
	\begin{theorem}
For every $n\geq 2$, $D_l(F_n)=D(F_n)= \lceil \dfrac{1+\sqrt{8n+1}}{2}\rceil .$  
\end{theorem}
	\proof It is sufficient to prove that $D_l(F_n)\leq D(F_n)$. For this purpose, we suppose that $L=\{L(v)\}_{v\in V(F_n)}$ is an arbitrary list assignment to $F_n$ in which $|L(v)|=D(F_n)$ for all $v\in V(F_n)$. We find a distinguishing labeling of $F_n$ such that the label of each vertex $v$ comes from $L(v)$. By Figure \ref{friend}, we set $L_i= L(v_i)$ for $1\leq i \leq 2n$.  The set of 2-element sets $\{a_i,b_i\}$ such that $a_i\in L_{2i-1}$, $b_i\in L_{2i}$ or $a_i\in L_{2i}$, $b_i\in L_{2i-1}$ is denoted by $\{L_{2i-1}, L_{2i}\}$ for every $1\leq i \leq n$. It is clear that $|\{L_{2i-1}, L_{2i}\}|\geq \frac{D(F_n)(D(F_n)-1)}{2}$. On the other hand  $D(F_n)= \lceil \dfrac{1+\sqrt{8n+1}}{2}\rceil$, and so $n \leq \frac{D(F_n)(D(F_n)-1)}{2}$. Hence for any $1\leq i \leq n$, we can assign an element of  $\{L_{2i-1}, L_{2i}\}$, say $\{a_i,b_i\}$, to the $\{v_{2i-1}, v_{2i}\}$ such that $\{a_i,b_i\}\neq \{a_j,b_j\}$ for every $i,j \in \{1,\ldots , n\}$ with $i\neq j$. Since  the central vertex $w$ is fixed under each automorphism, so we label the vertex $w$ with an arbitrary label in $L(w)$. This labeling is a distinguishing labeling of $F_n$. In fact, if $f$ is an automorphism of $F_n$ preserving the labeling, then $f$ maps the set $\{v_{2i-1}, v_{2i}\}$ to itself for every $1\leq i \leq n$, due to $\{a_i,b_i\}\neq \{a_j,b_j\}$ for every $i,j \in \{1,\ldots , n\}$ with $i\neq j$. Since the labels of $v_{2i-1}$ and $v_{2i}$ are distinct, so $f(v_{2i-1})=v_{2i-1}$ and $f(v_{2i})=v_{2i}$ for every $1\leq i \leq n$. Then $f$ is the identity automorphism of $F_n$ and the  labeling is distinguishing. So we have the result. \qed

\begin{figure}
	\begin{center}
		\hspace{.7cm}
		\includegraphics[width=0.5\textwidth]{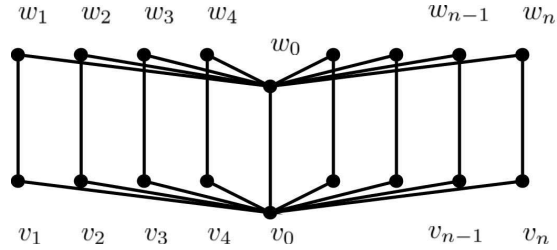}
				\caption{\label{book}  Book graph $B_n$.}
	\end{center}
\end{figure}

	\medskip
 The $n$-book graph $(n\geqslant 2)$ (Figure \ref{book}) is defined as the Cartesian product of $K_{1,n}$ and $P_2$, i.e., $K_{1,n}\square P_2$. We call every $C_4$ in the book graph $B_n$, a page of $B_n$. All pages in $B_n$ have a common side $v_0w_0$. We  computed the distinguishing number of $B_n$, and shall show that $D(B_n)=D_l(B_n)$. 

   \begin{theorem}{\rm \cite{Alikhani}}
	For every  $n\geq 2$, $D(B_n)= \lceil \sqrt{n}\rceil $.
\end{theorem}

   \begin{theorem}
	For every $n\geq 2$,  $D_l(B_n)=D(B_n)= \lceil \sqrt{n}\rceil$.
\end{theorem}
	\proof It is sufficient to prove that $D_l(B_n)\leq D(B_n)$. For this purpose, we suppose that $L=\{L(v)\}_{v\in V(B_n)}$ is an arbitrary list assignment to $B_n$ in which $|L(v)|=D(B_n)$ for all $v\in V(B_n)$. We find a distinguishing labeling of $B_n$ such that the label of each vertex $v$ comes from $L(v)$.  By Figure \ref{book}, the set of ordered pairs $(a_i,b_i)$ such that $a_i\in L(v_i)$ and $b_i\in L(w_i)$ is denoted by $(L(v_i), L(w_i))$ for every $1\leq i \leq n$. It is clear that $|(L(v_i), L(w_i))|\geq (D(B_n))^2$. On the other hand  $D(B_n)= \lceil \sqrt{n}\rceil$, and so $n \leq (D(B_n))^2$. Hence for any $1\leq i \leq n$, we can assign an element of  $(L(v_i), L(w_i))$, say $(a_i,b_i)$, to the $(v_{i}, w_{i})$ such that $(a_i,b_i)\neq (a_j,b_j)$ for every $i,j \in \{1,\ldots , n\}$ with $i\neq j$. Since  the central vertices $v_0$ and $w_0$ are fixed under each automorphism, so we label the vertices $v_0$ and $w_0$ with two different   labels in $L(v_0)$ and $L(w_0)$. This labeling is a distinguishing labeling of $B_n$. In fact, if $f$ is an automorphism of $B_n$ preserving the labeling, then $f(v_0)=v_0$ and $f(w_0) = w_0$. Thus $f$ maps  the set $\{(v_{i}, w_{i}): 1\leq i \leq n\}$ to itself. But since  $(a_i,b_i)\neq (a_j,b_j)$ for every $i,j \in \{1,\ldots , n\}$ with $i\neq j$, so $f$ maps  the set $(v_{i}, w_{i})$  to itself for every $1\leq i \leq n$. Since $f$ preserves the adjacency relation and $f(v_0)=v_0$ and $f(w_0) = w_0$, so $f(v_i)=v_i$ and $f(w_i) = w_i$ for every $1\leq i \leq n$.  Then $f$ is the identity automorphism of $B_n$ and the labeling is distinguishing. Therefore we have  the result. \qed

\section{Characterization of graph $G$ with $D(G)=D_l(G)$}

In this section we shall obtain a necessary and sufficient condition for a graph $G$ such that $D(G)=D_l(G)$. To do this, first we need to state some notation and results from 
set theory in subsection 3.1. In Subsection 3.2 we characterize graphs $G$ whose the  distinguishing number and the list distinguishing number are equal. 

\subsection{Some notations and results from set theory}

We begin this subsection with the following definition: 
\begin{definition}
Let $f:\{a_1,\ldots , a_n\}\rightarrow \{b_1,\ldots , b_d\}$, $d\leq n$, be a function. A $(m,d)$-related sequence to $f$ is $L^{(f)}=\{L_i\}_{i=1}^{n}$ such that $f(a_i)\in L_i$, $|L_i|=d$ and $L_i\subseteq \{b_1,\ldots , b_m\}$ for every $1\leq i \leq n$ where $m\geq d$. 
\end{definition}

It is clear that we have ${m-1 \choose d-1}^n$, $(m,d)$-related sequences to $f$. If the set of all these sequences is denoted by $\mathcal{L}^{\{a_1,\ldots , a_n\}}_{(m,d)}(f)$, then $|\mathcal{L}^{\{a_1,\ldots , a_n\}}_{(m,d)}(f)|={m-1 \choose d-1}^n$.

\begin{definition}\label{def2}
Let $f_i$, $1\leq i \leq t$, be functions of $\{a_1,\ldots , a_n\}$ into  $d$-subsets of  $\{b_1,\ldots , b_m\}$, where $d\leq n$ and $d\leq m$. The set of all possible different  $(m,d)$-related sequence to $f_1, \ldots , f_t$ is denoted by $B_{(m,d)}^{\{a_1,\ldots , a_n\}}(f_1,\ldots , f_t)$. In fact,  $B_{(m,d)}^{\{a_1,\ldots , a_n\}}(f_1,\ldots , f_t)=\bigcup_{i=1}^t \mathcal{L}^{\{a_1,\ldots , a_n\}}_{(m,d)}(f_i)$.
\end{definition}

It is clear that $|B_{(d,d)}^{\{a_1,\ldots , a_n\}}(f_1,\ldots , f_t)|=1$. To obtain a characterization of graphs $G$ with $D(G)=D_l(G)$, we shall know  the number of elements in $B^{\{a_1,\ldots , a_n\}}_{(m,d)}(f_1,\ldots , f_t)$ for any $m\geq d$. For computing this value we state  the following notation.

\medskip

\noindent \textbf{Notation.}  Let $f_i$, $1\leq i \leq t$, be functions of $\{a_1,\ldots , a_n\}$ into  $d$-subsets of  $\{b_1,\ldots , b_m\}$, where $d\leq n$ and $d\leq m$. Let $f_{j_1}, f_{j_2}, \ldots , f_{j_{i-1}}, f_{j_i}$ be $i$ different functions such that $1\leq j_1 < j_2 < \cdots < j_{i-1}< j_i\leq t$. For every $1\leq p \leq i$, we set
\begin{equation}\label{notation}
n^{(p)}_{\{j_1,j_2,\ldots , j_{i-1},j_i\}}=\big|\{a\in \{a_1,\ldots , a_n\}:~|\{f_{j_1}(a), f_{j_2}(a), \ldots , f_{j_{i-1}}(a), f_{j_i}(a)\}|=p \}\big|.
\end{equation}  
\begin{proposition}\label{3.3}
Let $f_1,\ldots , f_t$  be functions of $\{a_1,\ldots , a_n\}$ into  $d$-subsets of  $\{b_1,\ldots , b_m\}$, where $d\leq n$ and $d\leq m$. Then $|B_{(m,d)}^{\{a_1,\ldots , a_n\}}(f_1,\ldots , f_t)|=\sum_{i=1}^{t}S_i$, such that $S_1={m-1 \choose d-1}^n$, and for every $i\geq 2$ we have

\begin{align*}
S_i &= {m-1 \choose d-1}^n-\\
&\sum_{j=1}^{i-1}{m-2 \choose d-2}^{n^{(2)}_{\{j,i\}}}{m-1 \choose d-1}^{n^{(1)}_{\{j,i\}}}+\\
&\sum_{1\leq j < k  < i}{m-3 \choose d-3}^{n^{(3)}_{\{j,k,i\}}}{m-2 \choose d-2}^{n^{(2)}_{\{j,k,i\}}}{m-1 \choose d-1}^{n^{(1)}_{\{j,k,i\}}}-\\
&\sum_{1\leq j < k < h < i}{m-4 \choose d-4}^{n^{(4)}_{\{j,k,h,i\}}}{m-3 \choose d-3}^{n^{(3)}_{\{j,k,h,i\}}}{m-2 \choose d-2}^{n^{(2)}_{\{j,k,h,i\}}}{m-1 \choose d-1}^{n^{(1)}_{\{j,k,h,i\}}}+\\
& ~~~~~~~~~\vdots \\
&+ (-1)^{i+1}{m-i \choose d-i}^{n^{(i)}_{\{1,2,\ldots , i-1,i\}}}{m-(i-1) \choose d-(i-1)}^{n^{(i-1)}_{\{1,2,\ldots , i-1,i\}}} \cdots {m-1 \choose d-1}^{n^{(1)}_{\{1,2,\ldots , i-1,i\}}}.\\
\end{align*}
\end{proposition}
\proof
It is clear that the number of different sequences related to the function $f_1$ is $S_1= {m-1 \choose d-1}^n$. For the function $f_2$, there are  ${m-1 \choose d-1}^n$ related sequences, but some of them are exactly the same as related sequences to $f_1$, i.e., $\mathcal{L}^{\{a_1,\ldots , a_n\}}_{(m,d)}(f_1)\cap \mathcal{L}^{\{a_1,\ldots , a_n\}}_{(m,d)}(f_2)\neq \emptyset$. Using notation  \eqref{notation}, we conclude that $$|\mathcal{L}^{\{a_1,\ldots , a_n\}}_{(m,d)}(f_1)\cap \mathcal{L}^{\{a_1,\ldots , a_n\}}_{(m,d)}(f_2)|= {m-2 \choose d-2}^{n^{(2)}_{\{1,2\}}}{m-1 \choose d-1}^{n^{(1)}_{\{1,2\}}}.$$ 
Thus there exist $S_2$ new sequences related to $f_2$   by the inclusion-exclusion principle, where
$$S_2= {m-1 \choose d-1}^n-{m-2 \choose d-2}^{n^{(2)}_{\{1,2\}}}{m-1 \choose d-1}^{n^{(1)}_{\{1,2\}}}.$$
By a similar argument, for the function $f_3$, there are  ${m-1 \choose d-1}^n$ related sequences, but some of them are exactly the same as related sequences to $f_1$ and $f_2$. Using notation \eqref{notation}, we obtain that
\begin{itemize}
\item $|\mathcal{L}^{\{a_1,\ldots , a_n\}}_{(m,d)}(f_1)\cap \mathcal{L}^{\{a_1,\ldots , a_n\}}_{(m,d)}(f_3)|= {m-2 \choose d-2}^{n^{(2)}_{\{1,3\}}}{m-1 \choose d-1}^{n^{(1)}_{\{1,3\}}}$.
\item $|\mathcal{L}^{\{a_1,\ldots , a_n\}}_{(m,d)}(f_2)\cap \mathcal{L}^{\{a_1,\ldots , a_n\}}_{(m,d)}(f_3)|= {m-2 \choose d-2}^{n^{(2)}_{\{2,3\}}}{m-1 \choose d-1}^{n^{(1)}_{\{2,3\}}}$.
\item $|\mathcal{L}^{\{a_1,\ldots , a_n\}}_{(m,d)}(f_1) \cap \mathcal{L}^{\{a_1,\ldots , a_n\}}_{(m,d)}(f_2)\cap \mathcal{L}^{\{a_1,\ldots , a_n\}}_{(m,d)}(f_3)|= {m-3 \choose d-3}^{n^{(3)}_{\{1,2,3\}}}{m-2 \choose d-2}^{n^{(2)}_{\{1,2,3\}}}{m-1 \choose d-1}^{n^{(1)}_{\{1,2,3\}}}$.
\end{itemize} 

Hence there exist $S_3$ new sequences related to $f_3$ by the inclusion-exclusion principle, where
\begin{align*}
S_3={m-1 \choose d-1}^n- &
\left( {m-2 \choose d-2}^{n^{(2)}_{\{1,3\}}}{m-1 \choose d-1}^{n^{(1)}_{\{1,3\}}}+
{m-2 \choose d-2}^{n^{(2)}_{\{2,3\}}}{m-1 \choose d-1}^{n^{(1)}_{\{2,3\}}}\right)\\
&+{m-3 \choose d-3}^{n^{(3)}_{\{1,2,3\}}}{m-2 \choose d-2}^{n^{(2)}_{\{1,2,3\}}}{m-1 \choose d-1}^{n^{(1)}_{\{1,2,3\}}}.
\end{align*}
In general, for the function $f_i$, there are  ${m-1 \choose d-1}^n$ related sequences, but some of them are exactly the same as related sequences to $f_1, \ldots , f_{i-1}$. Using  notation \eqref{notation} and the inclusion-exclusion principle, we conclude that  there exist $S_i$ new sequences related to $f_i$ where
\begin{align*}
S_i &= {m-1 \choose d-1}^n-\\
&\sum_{j=1}^{i-1}{m-2 \choose d-2}^{n^{(2)}_{\{j,i\}}}{m-1 \choose d-1}^{n^{(1)}_{\{j,i\}}}+\\
&\sum_{1\leq j < k  < i}{m-3 \choose d-3}^{n^{(3)}_{\{j,k,i\}}}{m-2 \choose d-2}^{n^{(2)}_{\{j,k,i\}}}{m-1 \choose d-1}^{n^{(1)}_{\{j,k,i\}}}-\\
&\sum_{1\leq j < k < h < i}{m-4 \choose d-4}^{n^{(4)}_{\{j,k,h,i\}}}{m-3 \choose d-3}^{n^{(3)}_{\{j,k,h,i\}}}{m-2 \choose d-2}^{n^{(2)}_{\{j,k,h,i\}}}{m-1 \choose d-1}^{n^{(1)}_{\{j,k,h,i\}}}+\\
& ~~~~~~~~~\vdots \\
&+ (-1)^{i+1}{m-i \choose d-i}^{n^{(i)}_{\{1,2,\ldots , i-1,i\}}}{m-(i-1) \choose d-(i-1)}^{n^{(i-1)}_{\{1,2,\ldots , i-1,i\}}} \cdots {m-1 \choose d-1}^{n^{(1)}_{\{1,2,\ldots , i-1,i\}}}.\\
\end{align*}
Therefore, $|B_{(m,d)}^{\{a_1,\ldots , a_n\}}(f_1,\ldots , f_t)|=\sum_{i=1}^{t}S_i$.\qed

Here, we present another method for computing $|B_{(m,d)}^{\{a_1,\ldots , a_n\}}(f_1,\ldots , f_t)|$ by a recurrence relation.
\begin{proposition}\label{recurrence}
Let $f_1,\ldots , f_t$  be functions of $\{a_1,\ldots , a_n\}$ into  $d$-subsets of  $\{b_1,\ldots , b_m\}$, where $d\leq n$ and $d\leq m$. Then
$$|B^{\{a_1,\ldots , a_n\}}_{(m+1,d)}(f_1,\ldots , f_t)|= \sum_{i=0}^{n}\sum_{q=1}^{{n\choose i}}|B^{\{a_{q_1},\ldots , a_{q_i}\}}_{(m,d-1)}(f'_{q1},\ldots , f'_{qt})| |B^{\{a_{q_{i+1}},\ldots , a_{q_n}\}}_{(m,d)}(f''_{q1},\ldots , f''_{qt})|,$$
 where $f'_{q1},\ldots , f'_{qt}$ are the restriction of $f_1,\ldots , f_t$ to the $q$-th $i$-subset of $\{a_1,\ldots a_n\}$, say $\{a_{q_1},\ldots , a_{q_i}\}$, and $f''_{q1},\ldots , f''_{qt}$ are the restriction of $f_1,\ldots , f_t$ to the set $\{a_1,\ldots a_n\}\setminus \{a_{q_1},\ldots , a_{q_i}\}=\{a_{q_{i+1}},\ldots , a_{q_n}\}$.
\end{proposition}
\proof  Let $f_1,\ldots , f_t$  be functions of $\{a_1,\ldots , a_n\}$ into  $d$-subsets of  $\{b_1,\ldots , b_{m+1}\}$, where $d\leq n$ and $d\leq m+1$. Let $A_q=\{a_{q_1},\ldots , a_{q_i}\}$ be $q$-th  $i$-subset of $\{a_1,\ldots , a_n\}$ where $0\leq i \leq n$ and $1\leq q \leq {n \choose i}$ such that $\{a_1,\ldots , a_n\}\setminus A_q=\{a_{q_{i+1}},\ldots , a_{q_n}\}$.  If  $f'_{q1},\ldots , f'_{qt}$ are the restriction of $f_1,\ldots , f_t$ to the set $A_q$, then by Definition \ref{def2}, there exist $|B^{\{a_{q_1},\ldots , a_{q_i}\}}_{(m,d-1)}(f'_{q1},\ldots , f'_{qt})|$, $(m,d)$-related sequences $L'=\{L'_k\}_{k=q_1}^{q_i}$ with $|L'_k|=d$ and $L'_k\subseteq \{b_1,\ldots , b_{m+1}\}$ such that $b_{m+1}\in L'_k$ for every $k\in \{q_1,\ldots , q_i\}$.   If  $f''_{q1},\ldots , f''_{qt}$ are  the restriction of $f_1,\ldots , f_t$ to the set $\{a_{q_{i+1}},\ldots , a_{q_n}\}$, then  there exist $|B^{\{a_{q_{i+1}},\ldots , a_{q_n}\}}_{(m,d)}(f''_{q1},\ldots , f''_{qt})|$, $(m,d)$-related sequences $L''=\{L''_k\}_{k=q_{i+1}}^{q_n}$ with $|L''_k|=d$ and $ L''_k\subseteq \{b_1,\ldots , b_m\}$ for every $k\in \{q_{i+1},\ldots , q_n\}$.  Since $L=\{L_k\}_{k=1}^{n}$ where $L_k=L'_k$ for $k\in \{q_1,\ldots , q_i\}$, and $L_k= L''_k$ for $k\in \{q_{i+1}, \ldots , q_n\}$, is an element of   $B^{\{a_1,\ldots , a_n\}}_{(m+1,d)}(f_1,\ldots , f_t)$, and every element of $B^{\{a_1,\ldots , a_n\}}_{(m+1,d)}(f_1,\ldots , f_t)$ is obtained in such a way,  so we have the result by the rules of product and sum. \qed



\subsection{Graph $G$ with $D(G)=D_l(G)$ }

In this subsection, we present a necessary and sufficient condition for a graph $G$ with $D(G)=D_l(G)$. 

Let $G$ be a graph with $V(G)=\{a_1,\ldots , a_n\}$ and  $D(G)=d$. Let $L=\{L_i\}_{i=1}^n$ be an arbitrary sequence such that $|L_i|=d$ and $L_i\subseteq \{1, \ldots ,m\}$ for some $m\geq d$ and every $1\leq i \leq n$. If $L$ is a distinguishing $L$-labeling of $G$ then there exists a distinguishing labeling $C$ of vertices of $G$ such that $C(v_i)\in L_i$ for all $1\leq i \leq n$. On the other hand, for every distinguishing labeling $C$, we can construct   ${m-1 \choose d-1}^n$ sequences $L^{(C)}=\{L^{(C)}_i\}_{i=1}^n$ such that $C(v_i)\in L^{(C)}_i$,   $|L^{(C)}_i|=d$ and $L^{(C)}_i\subseteq \{1, \ldots ,m\}$ for every $1\leq i \leq n$. We call such sequences, the $(m,d)$-related sequences to $C$. If we denote the set of all related sequences to $C$ by $\mathcal{L}^{\{a_1,\ldots , a_n\}}_{(m,d)}(C)$, then $|\mathcal{L}^{\{a_1,\ldots , a_n\}}_{(m,d)}(C)|= {m-1 \choose d-1}^n$. Let $\mathcal{L}(G,m)$ be the set of all distinguishing labeling of $G$ with at most $m$ labels $\{1,\ldots , m\}$. Set $\mathcal{L}(G,m)=\{C_1, \ldots , C_{t_m}\}$. We suppose that $B_{(m,d)}^{\{a_1,\ldots , a_n\}}(C_1, \ldots , C_{t_m})$ is the set of all those sequences $L=\{L_i\}_{i=1}^n$ such that $|L_i|=d$ and $L_i\subseteq \{1, \ldots , m\}$ which are constructed using the distinguishing labelings in $\mathcal{L}(G,m)$, i.e., $B_{(m,d)}^{\{a_1,\ldots , a_n\}}(C_1, \ldots , C_{t_m})=\bigcup_{i=1}^{t_m}\mathcal{L}^{\{a_1,\ldots , a_n\}}_{(m,d)}(C_i)$. By these statements we have the following theorem:

\begin{theorem}
Let $G$ be a graph with $V(G)=\{a_1,\ldots , a_n\}$ and the distinguishing number $D(G)=d$. Let $\mathcal{L}(G,m)=\{C_1, \ldots , C_{t_m}\}$ be the set of all distinguishing labeling of $G$ with at most $m$ labels $\{1,\ldots , m\}$ where $m\geq d$. An arbitrary sequence $L=\{L_i\}_{i=1}^n$ with  $|L_i|=d$ and $L_i\subseteq \{1, \ldots ,m\}$ for every $1\leq i \leq n$, is a distinguishing $L$-labeling of $G$, if and only if  $L\in B_{(m,d)}^{\{a_1,\ldots , a_n\}}(C_1, \ldots , C_{t_m})$.
\end{theorem}

\medskip
If the set of all sequences  $L=\{L_i\}_{i=1}^n$ with  $|L_i|=d$ and $L_i\subseteq \{1, \ldots ,m\}$ for every some $m\geq d$ and $1\leq i \leq n$ is denoted by $A$, then $|A|={m\choose d}^n $. It is clear that $B_{(m,d)}^{\{a_1,\ldots , a_n\}}(C_1, \ldots , C_{t_m})\subseteq A$.  We note that we can compute $|B_{(m,d)}^{\{a_1,\ldots , a_n\}}(C_1, \ldots , C_{t_m})|$, by Propositions \ref{3.3} and \ref{recurrence}. The following theorem gives the characterization of graph $G$ with $D(G)=D_l(G)$. 
 
\begin{theorem}
Let $G$ be a graph with $V(G)=\{a_1,\ldots , a_n\}$ and the distinguishing number $D(G)=d$. Let $\mathcal{L}(G,m)=\{C_1, \ldots , C_{t_m}\}$ be the set of all distinguishing labeling of $G$ with at most $m$ labels $\{1,\ldots , m\}$ where $m\geq d$. Then 
\begin{itemize} 
\item $D_l(G)={\rm min}\{d:  |B_{(m,d)}^{\{a_1,\ldots , a_n\}}(C_1, \ldots , C_{t_m})|\geq {m\choose d}^n ~\text{for all}~m\geq d\}$.
\item   $D_l(G)=D(G)=d$ if and only if  $|B_{(m,d)}^{\{a_1,\ldots , a_n\}}(C_1, \ldots , C_{t_m})|\geq {m\choose d}^n $ for all $m \geq d$.
\end{itemize}
\end{theorem}

\section{List chromatic number and list distinguishing chromatic number }
In this  section, we show  that we can apply the results  of Section 3 for other list numbers, including list  chromatic number $\chi_l$, and list distinguishing chromatic number $\chi_{D_l}$. 
In 1979, Erd\"os, Rubin and Taylor \cite{erdos} introduced a beautiful new direction of graph colorings. We will say that a graph $G$ is $k$-choosable if for any assignment of a set ( or ``list") $L_v$ of $k$ colors to each vertex $v$ of $G$, it is possible to select a color $\lambda_v\in L_v$ for each $v$ so that $\lambda_u \neq  \lambda_v$ if $u$ and $v$ are adjacent. The list chromatic number $\chi_l(G)$ is defined to be the least $k$ such that $G$ is $k$-choosable. Clearly, $\chi_l(G)\geq \chi(G)$  for any $G$ (by taking $L_v=\{1,2, \ldots , \chi(G)\}$  for all $v$). However, the difference between $\chi_l(G)$ and $\chi(G)$ can be arbitrarily large. In general, $\chi_l$ can not be bounded in terms of $\chi$. 
In 2006 Collins and
Trenk \cite{Collins} introduced the \textit{distinguishing chromatic number}  $\chi_D(G)$ of a graph $G$, as  the minimum number of labels
in a distinguishing labeling of $G$ that is also a proper coloring. While not explicitly introduced in \cite{Ferrara}, it is natural to consider a list analogue of the distinguishing chromatic number.
Ferrara et al. \cite{Ferrara2} say that $G$ is \textit{proper distinguishing $L$-labeling} if there is a distinguishing labeling $f$ of $G$ chosen from the lists such that $f$ is also
a proper coloring of $G$. The \textit{list distinguishing chromatic number}  $\chi_{D_l} (G)$ of $G$ is the minimum integer $k$ such that $G$ is properly
$L$-distinguishable for any assignment $L$ of lists with $|L(v)| \geq k$ for all $v$.

Let $G$ be a graph with $V(G)=\{a_1,\ldots , a_n\}$ and the distinguishing chromatic number (resp. chromatic number) $\chi_D(G)=d$ (resp. $\chi(G)=d$). Let $L=\{L_i\}_{i=1}^n$ be an arbitrary sequence such that $|L_i|=d$ and $L_i\subseteq \{1, \ldots ,m\}$ for some $m\geq d$ and every $1\leq i \leq n$. If $L$ is a proper distinguishing  $L$-labeling (resp. proper $L$-labeling) of $G$ then there exists a proper distinguishing labeling (resp. proper labeling ) $C$ of vertices of $G$ such that $C(v_i)\in L_i$ for all $1\leq i \leq n$. On the other hand, for every proper distinguishing labeling (resp. proper labeling ) $C$, we can construct   ${m-1 \choose d-1}^n$ sequences $L^{(C)}=\{L^{(C)}_i\}_{i=1}^n$ such that $C(v_i)\in L^{(C)}_i$,   $|L^{(C)}_i|=d$ and $L^{(C)}_i\subseteq \{1, \ldots ,m\}$ for every $1\leq i \leq n$. We call such sequences, the $(m,d)$-related sequences to $C$. If we denote the set of all related sequences to $C$ by $\mathcal{L}^{\{a_1,\ldots , a_n\}}_{(m,d)}(C)$, then $|\mathcal{L}^{\{a_1,\ldots , a_n\}}_{(m,d)}(C)|= {m-1 \choose d-1}^n$. Let $\mathcal{L}(G,m)$ be the set of all proper distinguishing labeling (resp. proper labeling) of $G$ with at most $m$ labels $\{1,\ldots , m\}$. Set $\mathcal{L}(G,m)=\{C_1, \ldots , C_{t_m}\}$. We suppose that $B_{(m,d)}^{\{a_1,\ldots , a_n\}}(C_1, \ldots , C_{t_m})$ is the set of all those sequences $L=\{L_i\}_{i=1}^n$ such that $|L_i|=d$ and $L_i\subseteq \{1, \ldots , m\}$ which are constructed using the proper distinguishing labelings (resp. proper labelings) in $\mathcal{L}(G,m)$, i.e., $B_{(m,d)}^{\{a_1,\ldots , a_n\}}(C_1, \ldots , C_{t_m})=\bigcup_{i=1}^{t_m}\mathcal{L}^{\{a_1,\ldots , a_n\}}_{(m,d)}(C_i)$. Therefore we can conclude the following results by the same argument as Section 3:

\begin{theorem}
Let $G$ be a graph with $V(G)=\{a_1,\ldots , a_n\}$ and the distinguishing chromatic number  $\chi_D(G)=d$. Let $\mathcal{L}(G,m)=\{C_1, \ldots , C_{t_m}\}$ be the set of all proper distinguishing labeling  of $G$ with at most $m$ labels $\{1,\ldots , m\}$ where $m\geq d$. 
  If \linebreak $B_{(m,d)}^{\{a_1,\ldots , a_n\}}(C_1, \ldots , C_{t_m})$ is the set of all those sequences $L=\{L_i\}_{i=1}^n$ such that $|L_i|=d$ and $L_i\subseteq \{1, \ldots , m\}$ which are constructed using the proper distinguishing labelings  in $\mathcal{L}(G,m)$, then 
 \begin{itemize}
 \item An arbitrary sequence $L=\{L_i\}_{i=1}^n$ with  $|L_i|=d$ and $L_i\subseteq \{1, \ldots ,m\}$ for every $1\leq i \leq n$, is a proper distinguishing $L$-labeling  of $G$, if and only if  $L\in B_{(m,d)}^{\{a_1,\ldots , a_n\}}(C_1, \ldots , C_{t_m})$.
\item  $\chi_{D_l}(G)={\rm min}\{d:  |B_{(m,d)}^{\{a_1,\ldots , a_n\}}(C_1, \ldots , C_{t_m})|\geq {m\choose d}^n ~\text{for all}~m\geq d\}.$
 \item $\chi_{D_l}(G)=\chi_{D}(G)=d$ if and only if  $|B_{(m,d)}^{\{a_1,\ldots , a_n\}}(C_1, \ldots , C_{t_m})|\geq {m\choose d}^n $ for all $m \geq d$.
 \end{itemize}
 \end{theorem}
 
 We end this paper with the following theorem:
 
 \begin{theorem}
Let $G$ be a graph with $V(G)=\{a_1,\ldots , a_n\}$ and the chromatic number $\chi(G)=d$. Let $\mathcal{L}(G,m)=\{C_1, \ldots , C_{t_m}\}$ be the set of all  proper labeling of $G$ with at most $m$ labels $\{1,\ldots , m\}$ where $m\geq d$. 
 If $B_{(m,d)}^{\{a_1,\ldots , a_n\}}(C_1, \ldots , C_{t_m})$ is the set of all those sequences $L=\{L_i\}_{i=1}^n$ such that $|L_i|=d$ and $L_i\subseteq \{1, \ldots , m\}$ which are constructed using the proper labelings in $\mathcal{L}(G,m)$, then
\begin{itemize}
\item an arbitrary sequence $L=\{L_i\}_{i=1}^n$ with  $|L_i|=d$ and $L_i\subseteq \{1, \ldots ,m\}$ for every $1\leq i \leq n$, is a proper  $L$-labeling of $G$, if and only if  $L\in B_{(m,d)}^{\{a_1,\ldots , a_n\}}(C_1, \ldots , C_{t_m})$.
 \item $\chi_{l}(G)={\rm min}\{d:  |B_{(m,d)}^{\{a_1,\ldots , a_n\}}(C_1, \ldots , C_{t_m})|\geq {m\choose d}^n ~\text{for all}~m\geq d\}.$
\item $\chi_{l}(G)=\chi(G)=d$ if and only if  $|B_{(m,d)}^{\{a_1,\ldots , a_n\}}(C_1, \ldots , C_{t_m})|\geq {m\choose d}^n $ for all $m \geq d$.
\end{itemize}
\end{theorem}


\end{document}